\documentclass[11pt, a4paper]{article}
\usepackage{t1enc}
\usepackage[latin1]{inputenc}
\usepackage[english]{babel}
\usepackage{amsmath,amsthm}
\usepackage{mathrsfs}
\usepackage{amsfonts}
\usepackage{latexsym}
\usepackage[dvips]{graphicx}
\usepackage{float}
\usepackage[natural]{xcolor}
\usepackage{algorithm}
\usepackage{algorithmic}
\usepackage{enumerate,enumitem}
\usepackage{multirow}
\usepackage{xcolor}
\usepackage[colorlinks,linkcolor=blue]{hyperref}
\usepackage{lineno}
\usepackage{setspace}
\usepackage{fullpage}
\usepackage[title]{appendix}
\usepackage{todonotes}
\usepackage{amssymb}
\usepackage{array}

\usepackage{listings}
\usepackage{bbm}

\usepackage{caption}
\usepackage{subfigure}
\usepackage{url}
\usepackage{fullpage}
\usepackage{hyperref}
\usepackage[margin=2.3cm]{geometry}
\usepackage{enumitem,verbatim}

\newtheorem{theorem}{Theorem}[section]

\newtheorem{conjecture}[theorem]{Conjecture}

\theoremstyle{definition}
\newtheorem{problem}[theorem]{Problem}
\newtheorem{definition}[theorem]{Definition}

\usepackage{todonotes}
\title{Transversal Structures in Graph Systems: A Survey}

\author{
Wanting Sun \thanks{Data Science Institute, Shandong University, Jinan, 250100, China. Email: {\tt wtsun@sdu.edu.cn}. Supported by the China Postdoctoral Science Foundation (2023M742092).}
\and
Guanghui Wang \thanks{School of Mathematics, Shandong University, Jinan, 250100, China. Email: {\tt ghwang@sdu.edu.cn}. Supported by the National Key Research and Development Program (2023YFA1009603) and the Natural Science Foundation of China (12231018).
}
\and 
Lan Wei \thanks{School of Mathematics, Shandong University, Jinan, 250100, China. Email: {\tt lanwei@mail.sdu.edu.cn}.}
}
\date{}

\begin{document}
\maketitle
\begin{abstract}
    Given a system $\mathcal{G} =\{G_1,G_2,\dots,G_m\}$ of graphs/digraphs/hypergraphs on the common vertex set $V$ of size $n$, an $m$-edge graph/digraph/hypergraph $H$ on 
$V$ is \textit{transversal} in $\mathcal{G}$ if there exists a bijection $\phi :E(H)\rightarrow [m]$ such that 
$e \in E(G_{\phi(e)})$ for all $e\in E(H)$. In this survey, we consider extremal problems for transversal structures in graph systems. More precisely, we summarize some sufficient conditions 
that ensure the existence of transversal structures in graph/digraph/hypergraph systems, which generalize several classical theorems in extremal graph theory to transversal version. We also include a number of conjectures and open problems.
\end{abstract}
\maketitle

\section{Introduction}
\label{sec:1}
How {do} global parameters of a graph, such as its degree or chromatic number, influence its local substructures? 
Is there some sufficiently large average degree or chromatic number that forces a given substructure to appear? 
Questions of this type are one of the most natural problems in graph theory, and there are many profound and interesting results. These are collectively referred to as \textit{extremal graph theory}. It is a branch of combinatorial mathematics that studies the maximal or minimal properties of graphs under certain conditions. 

Many important problems in extremal graph theory can  be framed as subgraph containment problems, which ask for conditions on a graph $G$ that ensure it contains a copy of a general graph $H$. 
For example, Mantel's theorem states that every $n$-vertex graph with more than $\lfloor n^2/4\rfloor$ edges contains a triangle; Tur\'{a}n  \cite{Turan1941} generalized it  by determining the maximum number of edges in {an $n$-vertex} graph that forbids $K_r$, where $K_r$ denotes a complete graph on $r$ vertices;  
Dirac \cite{1952dirac} proved that every $n$-vertex ($n\geq 3$) graph  with  minimum degree at least 
$n/2$ contains a Hamilton cycle.  
For more results, the reader is referred to the book of Bollob\'{a}s  \cite{EGT}.


In this paper, we survey some recent developments on transversal problems in graph systems.
The notion of transversal appears in diverse branches of mathematics. A \textit{transversal} $X$ over a system $\mathcal{F}=\{F_1,\ldots,F_m\}$ of objects refers to  an object that intersects each $F_i$ exactly once. Here, the object can be any mathematical object, such as sets, spaces, set systems, and matroids, etc. Several classical theorems, which guarantee the existence of certain objects, have been reiterated in the context of transversals, demonstrating the possibility and impossibility of obtaining such an object as a transversal over a certain system. To name a few, transversals are extensively studied for Carath\'eodory's theorem \cite{1982Barany}, Helly's theorem \cite{2005Kalai_Helly},  Erd\H{o}s--Ko--Rado theorem \cite{2017AhaEKR}, and Rota's basis conjecture \cite{1994HuangLatin,2020PokrovskiyRota}, etc. 

Building on this trend, there have been numerous studies dedicated to generalizing classical results in extremal graph theory to transversal {versions}.
In addition to being interesting in their own right, they often provide a strengthening of the original results, which also reflects the robustness of graph properties. 
A transversal in a graph system  is also referred to as a ``rainbow subgraph'' in much literature, and it was formally
introduced in the recent work of Joos and Kim  \cite{2021jooskim}.  For a positive integer $n$, we write $[n] := \{1, 2,\dots, n\}$. 


\begin{definition}
     For a given system $\mathcal{G}=\{G_1,\ldots,G_m\}$ of graphs/digraphs/hypergraphs with the same vertex set $V$, an $m$-edge graph/digraph/hypergraph $H$ with vertices in $V$ is a \textit{$\mathcal{G}$-transversal} (or a transversal inside $\mathcal{G}$) if there exists a bijection $\phi:E(H)\to[m]$ such that $e\in E(G_{\phi(e)})$ for all $e\in E(H).$
\end{definition}

By interpreting each $G_i$ as the set of edges colored with the color $i$, 
the function $\phi$ in the above definition is referred to as a \textit{rainbow coloring} and $H$ is a \textit{rainbow subgraph} inside  $\mathcal{G}$.
We shall emphasize that in what follows, when we say that $\mathcal{G}$ is a graph system with order $n$, then  all graphs in $\mathcal{G}$ are not necessarily distinct and they have the same vertex set $V$ of size $n$. Many classical results in extremal graph theory have been extended to transversal
settings, exhibiting interesting phenomena.

Perhaps the most famous transversals are those related to Latin squares, which were studied by Euler. In 1782, Euler \cite{1782Euler} considered a Latin square of order $n$, which is an $n\times n$ array filled with symbols $1,\ldots,n$, where each symbol appears exactly once in each row and column, see Figure~(a). 
A \textit{partial transversal} of a Latin square is a set of cells such that every pair of cells shares no row, column or symbol, in particular a \textit{transversal} is a partial transversal with $n$ cells. 
The famous Ryser--Brualdi--Stein conjecture \cite{Brualdi,1967Ryser,Stein} is one of the most important problems on transversals in Latin squares.

\begin{figure}[htbp]
  \centering

\subfigure[Latin square of order $4$.]{
    \includegraphics[width=40mm]{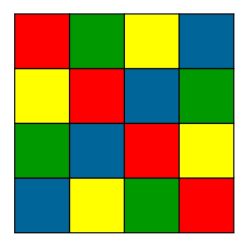}}\ \ \ \ \ \ \ \ \ \      
\subfigure[Edge-colored bipartite graph.]{
    \includegraphics[width=50mm]{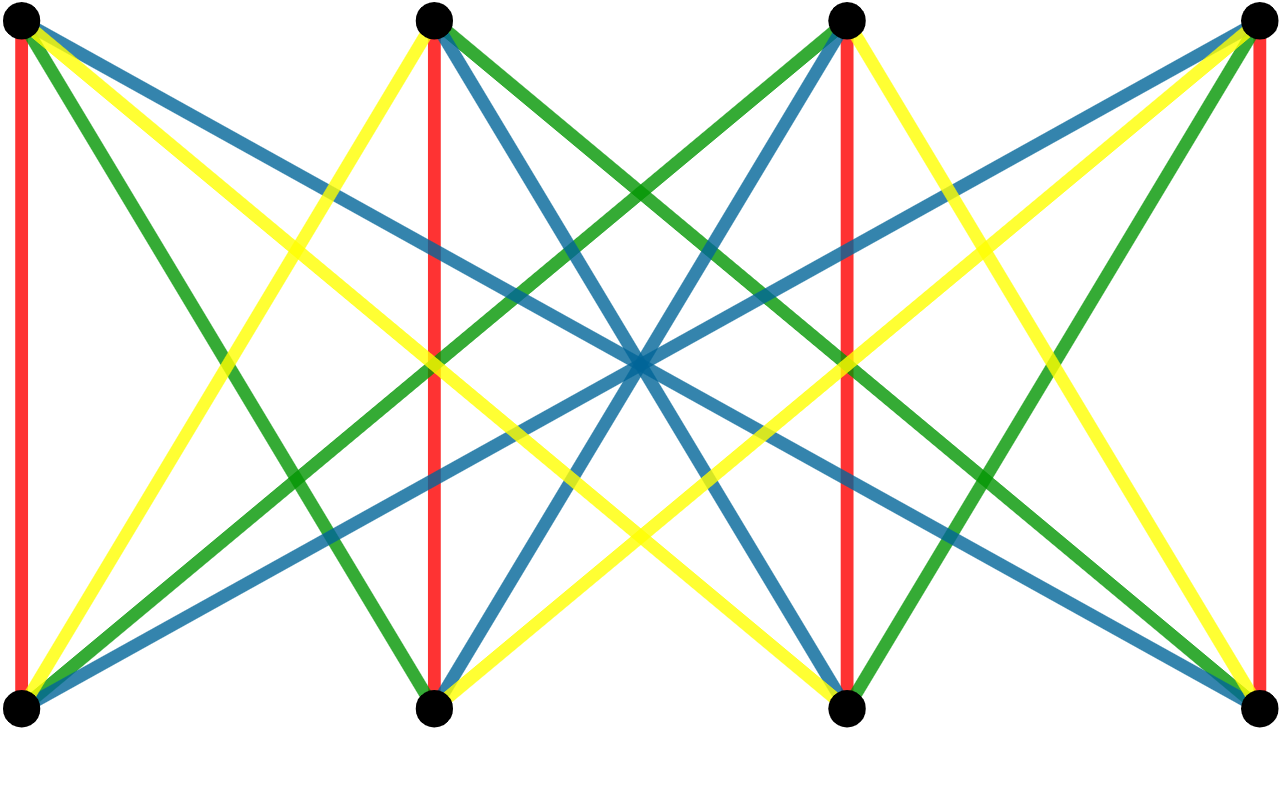}}

\end{figure}

\begin{conjecture}[Ryser--Brualdi--Stein conjecture, \cite{Brualdi,1967Ryser,Stein}]
 Every Latin square
of order $n$ has a partial transversal with $n-1$ cells, and a transversal with $n$ cells if $n$ is odd.
\end{conjecture}

In 2023, Montgomery \cite{montgomery2023} introduced the first techniques to identify and exploit the possible algebraic properties behind the entries in a Latin square, and proved the following. 
\begin{theorem}[\cite{montgomery2023}]
    There is some $n_0 \in \mathbb{N}$ such that every Latin square of order $n \geq n_0$ contains a partial transversal with $n - 1$ cells.
\end{theorem}
For more relevant results of transversals in Latin squares, see the recent survey by Montgomery \cite{2024Monlatin},  which collected results from the past decade. 

Considering the rows and columns of the Latin square as a complete  bipartite graph $K_{n,n}$, 
where each symbol represents a color and each cell represents a colored edge, the Latin square naturally corresponds to a properly edge-colored $K_{n,n}$, see Figure (b). The edge-colored $K_{n,n}$ can be viewed as a set of graphs $\mathcal{G}=\{G_1,\ldots, G_n\}$, where each $G_i$ consists of edges with color $i$, then a rainbow matching of $\mathcal{G}$ corresponds to a transversal of the Latin square. We refer the reader to see the nice survey \cite{Wanless2011}. 
These studies have inspired research on the existence of transversal structures in a system of graphs. 
The central question in the study of transversals can be generalized and depicted as follows.


\begin{problem}[\cite{2021jooskim}]
    Let $H$ be a graph/digraph/hypergraph with $m$ edges, and let $\mathcal{G}=\{G_1,G_2,\dots,\linebreak G_m\}$ be a system of graphs/digraphs/hypergraphs. Which property imposed on $\mathcal{G}$ will ensure a transversal copy of $H$?
\end{problem}

Since we do not require the objects in the family $\mathcal{G}$ to be distinct, a necessary condition for a positive answer to the
above problem is that every $n$-vertex graph/digraph/hypergraph satisfying the property contains a copy of $H$. {In this survey}, it will be shown by some examples that this necessary condition is not sufficient in some cases. 


In this paper, we aim to survey known results for transversal problems in graph/{random graph}/digraph/hypergraph  systems. Throughout the survey, we
include a number of conjectures and open problems.
It is of course impossible to cover all the known results for graph systems in a single survey; thus, the choice of results we present is inevitably subjective. Nevertheless, we hope to describe enough results from this fascinating area to engage many researchers in extremal and probabilistic combinatorics, and motivate further study of this subject.

\section{ Tur\'{a}n-type problem in graph systems}

One of the central topics in extremal combinatorics is to determine the maximum number of edges of a graph on $n$ vertices that does not contain a copy of a given graph $H$ as a subgraph, such a maximum number is denoted by $\mathrm{ex}(n,H)$. 
The most notable result is Tur\'{a}n's theorem \cite{Turan1941}, which  states that $\mathrm{ex}(n, K_{r+1}) = t_{r}(n)$, where $t_r(n)$ is the number of edges of a balanced complete $r$-partite graph on $n$ vertices (denoted by $T_r(n)$). Erd\H{o}s, Stone and Simonovits \cite{ESS,erdos-stone} found {a connection} between $\mathrm{ex}(n,H)$ and the chromatic number $\chi(H)$ of $H$. 


In addition to the Tur\'{a}n problem itself, there are also many variants involving colorings. For example, Diwan and Mubayi \cite{2006Diwan} considered a generalization of Tur\'an's theorem in a graph system containing two graphs. 
It is of course possible to consider similar {questions} with more than two graphs. In this section, we survey results on the Tur\'{a}n-type problems in graph systems. In other words, we
are to determine the ``largest'' $n$-vertex graph systems that do not contain a copy of $H$ as a rainbow subgraph.  

Define the $k$-\textit{rainbow extremal number} $\mathrm{ex}_k^*(n,H)$ of a graph  $H$ to be the maximum $\min_{i\in [k]}|E(G_i)|$ among all rainbow $H$-free graph systems $\mathcal{G}=\{G_1,\ldots,G_k\}$ with order $n$. 
The simplest case here is a rainbow triangle. 
Diwan and Mubayi \cite{2006Diwan} proposed a problem: determine the value of  $\mathrm{ex}_3^*(n,K_3)$. 
Magnant \cite{2015MagnantTra} answered this problem under the assumption that the multigraphs (i.e., graphs formed by taking the union of the edge sets $E(G_i)$) in question contain at least one edge between every pair of vertices.
In 2020, Aharoni, DeVos, de la Maza, Montejano and \v{S}\'amal \cite{2019Aharoni} solved the above problem completely. 
\begin{theorem}[\cite{2019Aharoni}]\label{thmK3}
   $\mathrm{ex}_3^*(n,K_3)=\lfloor\frac{26-2\sqrt{7}}{81}\rfloor n^2$.
\end{theorem}

 Interestingly, the constant $\frac{26-2\sqrt{7}}{81}$ is an irrational number and larger than $\frac{1}{4}$.
Therefore,
compared to Mantel's theorem, the transversal triangle problem requires a stronger size 
condition.
This signifies the difference
between the transversal setting and the classical setting, and renders the rainbow extremal number more intriguing. Aharoni, DeVos, de la Maza, Montejano and \v{S}\'amal \cite{2019Aharoni} proposed a potential interesting direction to proceed.
\begin{problem}[\cite{2019Aharoni}]\label{triangle}
    For what real numbers $\alpha_1,\alpha_2,\alpha_3>0$, { is it true that  every graph system $\{G_1,G_2,G_3\}$ satisfying $|E(G_i)|>\alpha_i n^2$ 
    must have a transversal  triangle?}
\end{problem}
Falgas-Ravry, Markstr\"om and R\"aty \cite[Theorem 1.13]{2024Falgas-RavryTra} solved this problem when $n$ is sufficiently large. 
Recall that Tur\'{a}n's theorem \cite{Turan1941} generalizes Mantel's theorem by determining the maximum number of edges in $n$-vertex $K_r$-free ($r\geq 3$) graphs. 
Analogously, one may consider the
following.
\begin{problem}[\cite{2019Aharoni}]
Determine the value of $\mathrm{ex}_{\binom{r}{2}}^*(n,K_r)$ for $r\geq 4$.
\end{problem}

In general cases, the problem becomes more complicated and there are only few known results related to the $k$-rainbow extremal number $\mathrm{ex}_k^*(n,H)$. 
As a corollary of Keevash, Saks, Sudakov and Verstra\"{e}te \cite{2004KeevashMulticolour},
we have $\mathrm{ex}_k^*(n,K_3)\leq n^2/4 +o(n^2)$ when $k\geq 4$. Frankl \cite{2022FranklRTra} gave a new proof of this result. Babi\'{n}ski and Grzesik \cite{2024BabinskiP3} determined the value of $\mathrm{ex}_k^*(n,P_4)$  up to {an} additive $o(n^2)$ error term for every $k\geq 3$, where $P_4$ denotes the path on $4$ vertices.

The \textit{$k$-rainbow Tur\'{a}n density} of a graph $H$ is defined as follows:
$$
 \pi_k^*(H)=\lim_{n\rightarrow \infty}\frac{\mathrm{ex}_k^*(n,H)}{\binom{n}{2}}. 
$$
Im, Kim, Lee and Seo \cite{Im2023} showed that the limit exists for every graph $H$. They \cite{Im2023} conjectured that the rainbow Tur\'{a}n density is always an  algebraic number, and they proved the conjecture for all trees. 
Furthermore, they \cite{Im2023}  gave an answer to an analogue of Problem \ref{triangle} for trees.
\begin{theorem}[\cite{Im2023}]
    Let $\alpha_1,\ldots,\alpha_k$ be nonnegative reals such that $\sum_{i=1}^{k}(1-\sqrt{\alpha_i})<1$. For every $\epsilon>0$, there
exists an integer $n_0$ such that for every $n\geq n_0$ and every $k$-edge tree $T$, every graph system $\mathcal{G}=\{G_1,\ldots,G_k\}$ of order $n$ with $|E(G_i)|\geq (\alpha_i+\epsilon){\binom{n}{2}}$ contains a transversal copy of $T$.
\end{theorem}

Except for the rainbow Tur\'{a}n number of a graph, depending on the contexts, different measures for the ``size'' of graph systems have been introduced, for example,  the sum $\sum_{i}|E(G_i)|$ \cite{2024Liucritical,2004KeevashMulticolour}, the product $\Pi_{i}|E(G_i)|$ \cite{2022FranklRTra,2022Frankldis}  and other measures. 

The \textit{$k$-color  Tur\'{a}n number}, denoted by $\mathrm{ex}_k(n,H)$, is the maximum value of $\sum_{i\in [k]}|E(G_i)|$ among all graph systems $\mathcal{G}=\{G_1,\ldots,G_k\}$ of order $n$ and without rainbow copy of $H$. We call such a graph system $\mathcal{G}$ satisfying $\sum_{i\in [k]}|E(G_i)|=\mathrm{ex}_k(n,H)$ the \textit{$k$-color extremal graph system of $H$}. 
To determine $\mathrm{ex}_k(n,H)$, there are two natural lower bounds that arise from the following constructions. One is that
each $G_i$ is the $n$-vertex $H$-free graph with $\mathrm{ex}(n,H)$ edges;
the other one is that there are $|E(H)|-1$ copies of the complete graph $K_n$ and all the 
other graphs are empty. 

A graph is called \textit{$r$-color-critical} if it has chromatic number $r$, and it has an edge whose removal reduces the chromatic number to $r-1.$  In 2004, Keevash, Saks, Sudakov  and Verstra\"{e}te \cite{2004KeevashMulticolour} proposed the following conjecture.

\begin{conjecture}[\cite{2004KeevashMulticolour}]\label{2004keevashMul}
    Suppose $r\geq 3$ and $k\geq h$. Let $H$ be an $r$-color-critical graph with $h$ edges. Then, there exists an $n_{0}=n_{0}(H)>0$ such that for all $n\geq n_0$, the $k$-color extremal graph system of $H$ either consists of exactly $h-1$ nonempty graphs where each of them is a copy of $K_n$  or consists of $k$ graphs where all of them are identical copies of $T_{r-1}(n)$. In particular, 
$$\mathrm{ex}_{k}( n, H)=
\begin{cases}(h-1)\binom{n}{2},&\ \text{if}\ h\leq k<\frac{r-1}{r-2}(h-1),\\
k\cdot t_{r-1}(n),&\ \text{if}\ k\geq\frac{r-1}{r-2}(h-1).
\end{cases}$$
\end{conjecture}
They provided evidence for this conjecture by
proving it holds when $H$ is a complete graph or $r=3$. Recently, Chakraborti, Kim, Lee, Liu and Seo \cite{2024Liucritical} showed that the conjecture holds for every $4$-color-critical graph $H$ and almost all $r$-color-critical
graphs $H$ with $r\geq 5$. Furthermore, they \cite{2024Liucritical} proved that if $H$ is not $r$-color-critical for some $r \geq 3$ and $n$ is
sufficiently large, then neither of  the constructions in Conjecture \ref{2004keevashMul} is extremal. Li, Ma and Zheng \cite{2024LicriticalConj} proved Conjecture \ref{2004keevashMul} holds  when $k\geq2(r-1)(h-1)/r$, significantly improving earlier results. Furthermore, their proof combined the stability argument of Chakraborti, Kim, Lee, Liu and Seo \cite{2024Liucritical} with a novel graph packing technique for embedding multigraphs.
\begin{theorem}[\cite{2024LicriticalConj}]
    Let $r\geq 5$ and $H$ be an $r$-color-critical graph with $h$ edges. If $n$ is sufficiently large and $k\geq\frac{2(r-1)(h-1)}{r}$, then $\mathrm{ex}_{k}( n, H) = k\cdot t_{r- 1}( n)$, and the unique $n$-vertex $k$-color extremal graph system of $H$ consists of $k$ colors all of which are identical copies of $T_{r-1}(n).$
\end{theorem}

There is also some research focusing on maximizing the product  $\Pi_{i}|E(G_i)|$ of the number of edges of a graph system without a rainbow copy of $H$. Frankl \cite{2022FranklRTra} conjectured that a graph system $\{G_1,G_2,G_3\}$ on a common vertex set of order $n$ with $\Pi_{i\in [3]}|E(G_i)|> \lfloor{n^2/4}\rfloor^3$ contains a transversal triangle.  This conjecture was disproved by Frankl, Gy\H{o}ri, He, Lv, Salia, Tompkins, Varga and Zhu \cite{2022Frankldis}, and they further 
{provided the following conjecture.}
\begin{conjecture}[\cite{2022Frankldis}]\label{2022Frankldis}
Let $\mathcal{G}=\{G_1,G_2,G_3\}$ be 
a graph system of order $n$. If $\mathcal{G}$ contains no transversal triangle, then $$\Pi_{i\in [3]}|E(G_i)|\leq\gamma n^{6}(1+o(1)),$$
where $\gamma$ is the maximum of
$\frac{x^2}{2}\left(\frac{x^2}{2}+\frac{(1-x)^2}{2}\right)\left(x(1-x)+\frac{(1-x)^2}{2}\right)$ on $[0,1]$.
\end{conjecture}

{Recall that Falgas-Ravry, Markstr\"om and R\"aty \cite{2024Falgas-RavryTra} solved Problem \ref{triangle}  for sufficiently large $n$. As a direct consequence, Conjecture \ref{2022Frankldis} holds when $n$ is sufficiently large.} 



\section{Dirac-type problem
in graph systems}

Let $G = (V,E)$ be a graph. For a vertex $v\in V$, denote by $N_G(v)$ the neighborhood of $v$ 
and let $d_G(v)=|N_G(v)|$ be the degree of $v$. 
The subscripts are omitted when $G$ is clear from the context. Let $\delta(G)$ and $\Delta(G)$ be the minimum degree and maximum degree of $G$ respectively. Let $\mathcal{G} =\{G_1,G_2,\dots,G_m\}$ be a graph system. Denote by  $\delta(\mathcal{G})=\min\{\delta(G_i):i\in [m]\}$ the minimum degree of $\mathcal{G}$.

In this section,  we  investigate the minimum degree
condition that {ensures} the existence of  transversal spanning  structures in graph systems. When one is interested in embedding a spanning graph, it makes sense to consider minimum degree conditions rather than size conditions, as, for example, a graph could be very dense but contain an isolated vertex, and therefore not contain a copy of any given spanning graph with minimum degree at least $1$.

Given the difficulty of obtaining exact results in {(rainbow)} graph embedding problems, many results have focused on determining the minimum $d$ such that for all $\epsilon>0$ and sufficiently large integer $n$, each graph (system) of order $n$ with minimum degree at least $(d+\epsilon)n$ contains a {(rainbow)} copy of $H$. Note that  when we consider this problem in graph systems, $d$ is at least the value required to guarantee a copy of $H$ in a single graph, since all graphs in the system could be identical.


\begin{definition}\label{defthreshold}
Let $\mathbf{F}$ be an infinite family of graphs. Denote by $\delta_{\mathbf{F}}$ (if it exists) the smallest real number $\delta$ such that for all $\alpha>0$ and for all but finitely many $F\in \mathbf{F}$ the following holds. Let $n = |V (F)|$ and $H$ be any $n$-vertex graph with $\delta(H)\geq(\delta+\alpha)n$. Then $H$ contains a copy of $F$.

Denote by $\delta_{\mathbf{F}}^{rb}$ (if it exists) the smallest real number $\delta$ such that for all $\alpha>0$ and for all but finitely many $F\in \mathbf{F}$ the following holds. Let $n=|V(F)|$ and $m=|E(F)|$, and let $\mathcal{H}=\{H_1,\ldots, H_{m}\}$ be any graph system of order $n$ with $\delta(\mathcal{H})\geq (\delta+\alpha)n$. Then $\mathcal{H}$ contains a transversal copy of $F$.



\end{definition}


If the above two values exist, then
it is clear that $\delta_{\mathbf{F}}^\mathrm{rb}\geq \delta_{\mathbf{F}}.$ {Moreover}, Montgomery, M\"{u}yesser and Pehova \cite{2022Montgomery} observed that $\delta_{\mathbf{F}}^{rb}$ can be much larger than $\delta_{\mathbf{F}}$. Set $\mathbf{F}=\{k\times (K_{2,3}\cup C_4):k\in \mathbb{N}\}$, where $k\times G$ denotes the graph obtained by taking $k$ vertex-disjoint copies of $G$. K\"{u}hn and Osthus \cite{KuhnOsthus} showed that $\delta_{\mathbf{F}}=4/9$. Consider the graph system $\mathcal{H}=\{H_1,\ldots,H_m\}$ on $V$ obtained in the following way. Partition $V$ into two almost equal vertex subsets, say $A$ and $B$, and for each $i\in [m-1]$,  assume that $H_i$ is the disjoint union of a clique on $A$ and a clique on $B$. Suppose that $H_m$ is a  complete bipartite graph between $A$ and $B$. Observe that if $\mathcal{H}$ contains a transversal copy of some $F\in \mathbf{F}$, the edge of $K_{2,3}$ or $C_4$ that belongs to $H_m$ would be a bridge of $F$. However, neither $K_{2,3}$ nor $C_4$ contains a bridge. Hence, $\delta_{\mathbf{F}}^{rb}\geq \delta(\mathcal{H})\geq 1/2$.



On the other hand, there are many natural instances with $\delta_{\mathbf{F}}^\mathrm{rb}= \delta_{\mathbf{F}}$, for example, { the family of $n$-cycles}
(see \cite{2021jooskim}).  When this equality holds, we 
{say} that the corresponding $\mathbf{F}$ is \textit{color-blind}. 


\subsection{Hamilton path and Hamilton cycle}

A \textit{Hamilton path} (\textit{cycle}) in a graph is a path (cycle) that visits every vertex exactly once. 
Hamiltonicity is one of the most well-studied properties in graph theory.
Dirac \cite{1952dirac} showed that every graph {on $n\geq 3$} vertices with minimum degree at least $n/2$ contains a Hamilton cycle. {An $n$-vertex graph is a \textit{Dirac graph} if its minimum degree is at least $n/2$.} Later on, Ore \cite{1960Ore} relaxed this condition and proved that every $n$-vertex graph  $G$ with $\sigma_2(G)\geq n$ contains a Hamilton cycle, where {$n\geq 3$} and $\sigma_2(G)=\min\{d(x)+d(y):xy\notin E(G)\}$. 
 Moon and Moser \cite{Moon} established the bipartite {analogue} of Dirac's theorem.

In 2020, Aharoni \cite{2019Aharoni} conjectured that Dirac's theorem can be extended to a transversal version. Cheng, Wang and Zhao \cite{2021ChengWangZhao} proved an asymptotic version of
this conjecture; Joos and Kim  \cite{2021jooskim} confirmed the conjecture completely. 
 \begin{theorem}[\cite{2021jooskim}]\label{Hamilton-cycle}
     Let $n\in \mathbb{N}$  and $n\geq 3$. Suppose that $\mathcal{G} = \{ G_1, \ldots , G_n\}$ is a graph system of order $n$ and  $\delta(\mathcal{G})\geq \frac{n}{2}$. Then $\mathcal{G}$ contains a transversal Hamilton cycle.
 \end{theorem}

 Bradshaw \cite{2021Bradshaw} gave a degree condition to guarantee the existence of transversal Hamilton cycles in bipartite graph systems, which generalizes the result of Moon and Moser \cite{Moon}. 
 
The \textit{distance} of two distinct vertices in a graph is the length of a shortest path connecting them. {The \textit{$k$-th power} $G^k$ of a graph $G$ is the graph obtained
from $G$ by joining every pair of vertices with distance at most $k$ in $G$. P\'osa (for $r=3$, \cite{Posa}) and Seymour (for $r>3$, \cite{seymour}) conjectured that for all $r\geq3$ and all $n\geq r$, every $n$-vertex graph $G$ with $\delta(G)\geq {(r-1)}n/r$ must contain the $(r-1)$-th power of a Hamilton cycle. This conjecture has been confirmed by Koml\'{o}s, S\'{a}rk\"{o}zy, and Szemer\'{e}di \cite{1996KomlosPosaconj,Komlos-conj} {for sufficiently large $n$}. 
Recently, Gupta, Hamann, M\"{u}yesser, Parczyk and Sgueglia \cite{2023GuptaPowerH} determined the asymptotically tight bound on $\delta(\mathcal{G})$ for the existence of transversal powers of Hamilton cycles, and they proposed a problem for the exact version of transversal P\'osa--Seymour conjecture. 
\begin{problem}[\cite{2023GuptaPowerH}]
    Let $\mathcal{G}=\{G_1,\ldots,G_{(r-1)n}\}$ be a graph system of order $n$. If $\delta(\mathcal{G})\geq (1-\frac{1}{r})n$,  is {this} already sufficient for a transversal copy of the $(r-1)$-th power of a
Hamilton cycle?
\end{problem}
}

Bowtell, Morris, Pehova and Staden \cite{2023universality} proposed a further generalization of the transversal problem, where one seeks a copy of $H$ with a given coloring.
\begin{problem}[\cite{2023universality}]\label{pattern}
    Given a graph system $\mathcal{G}=\{G_1,\ldots,G_m\}$ {and a graph $H$}, what conditions imposed on  $\mathcal{G}$ will  guarantee that for every coloring $\chi: E(H)\rightarrow \{1,\ldots,m\}$, there is a copy of $H$ such that the image of $e$ is an edge of $G_{\chi(e)}$ for all $e\in E(H)$?
\end{problem}
When $H$ is a Hamilton cycle, Bowtell, Morris, Pehova and Staden \cite{2023universality} gave an asymptotically optimal minimum degree condition for Problem~\ref{pattern}. 

\begin{theorem}[\cite{2023universality}]\label{patterncycle}
     For every $\alpha > 0$ there exists $n_{0}$ such that for every $n\geq n_{0}$  the following holds. Let $m\in \mathbb{N}$  and let $\mathcal{G}= \{ G_1, \ldots , G_m\}$  be a graph system of order $n$. If $\delta ( \mathcal{G} ) \geq ( \frac{1}{2}+ \alpha)n$, then for every $\chi : [n] \to [m]$, there is a Hamilton cycle $C$ such that the image of $e$ is an edge of $G_{\chi(e)}$ for all $e\in E(C)$. 
\end{theorem}

Determining the exact minimum degree condition to guarantee all Hamilton cycle patterns remains an
intriguing open problem. Bowtell, Morris, Pehova and Staden  \cite{2023universality} provided 
a construction demonstrating that $n/2+1$ is a necessary lower bound and {conjectured this is} an exact threshold. In light of this, it seems more pertinent to either obtain results that improve this $o(n)$ term to an explicit and relatively small constant, or find counterexamples to this being possible. 


It is known that one can measure the robustness of Dirac graphs
 with respect to Hamiltonicity by counting  the number of distinct Hamilton cycles that a Dirac graph must contain. For more details about the robustness, we refer the reader to see the nice survey of Sudakov \cite{SudakovRobust}.  This method can be extended to measure the robustness of Dirac graph systems. 
 
Bradshaw, Halasz and Stacho \cite{Brad2022} showed that given a system of $n$ Dirac graphs with order $n$, there are at least $(cn/e)^{cn}$ different transversal Hamilton cycles for some constant $c\geq 1/68$. Theorem \ref{patterncycle} recovers a stronger lower bound, that there exist at least $n!\sim \sqrt{2\pi n}(n/e)^n$ different transversal Hamilton cycles in a graph system  $\mathcal{G}$ consisting of $n$ {asymptotically Dirac
graphs (i.e., $\delta(\mathcal{G})\geq (1/2+\alpha)n$ for $\alpha>0$)}. Anastos and Chakraborti \cite{2023Anastosrobust} improved this by showing that there exists a constant $C>0$ such that there are $(Cn)^{2n}$ different transversal Hamilton cycles in a graph  system consisting of $n$ Dirac graphs, which is tight up to the choice of $C$. They \cite{2023Anastosrobust} also showed that every graph system consisting of $n$ Dirac graphs with order 
$n$ spans $n/2$ pairwise edge-disjoint transversal Hamilton cycles, and this is best possible.

Cheng and Staden \cite{2024cheng} established a stability result for transversal Hamilton cycles. 
 Given integers $a,b\geq 0$, let $\mathcal{H}_a^b$ be a graph
system obtained by taking $a$ copies of $K_{\frac n2}\cup K_{\frac n2}$ and $b$ copies of $K_{{\frac n2},{\frac n2}}$ where
they are defined on the same equitable partition $A\cup B$. 
For two {disjoint} vertex subsets $X,\,Y \subseteq V$, 
let $G[X]$ denote the subgraph of 
$G$ induced by $X$, let $G[X,Y]$ be the subgraph of $G$ with vertex set $X\cup Y$ and edge set $\{xy:x\in X,\,y\in Y\ \text{and}\ xy\in E(G)\}$. 
We say that a graph system $\mathcal{J}$ on the vertex set 
$V$ of order $n$ is \textit{half-split} if there is $A\subseteq V$
with $|A|=\lfloor n/2\rfloor+1$ such that $J[A]=\emptyset$ and $J[A,V\setminus A]$ is complete for all $J\in \mathcal{J}$. Given $\kappa>0$, we say that two graph systems $\mathcal{G}$ and $\mathcal{G}'$ on a common vertex set of order $n$ are \textit{$\kappa$-close} if $\mathcal{G}'$ can be obtained from $\mathcal{G}$ by adding or deleting at most $\kappa n^3$ edges.   

 \begin{theorem}[\cite{2024cheng}]\label{stability}
     For all $\kappa>0$, there exist $\alpha>0$ and $n_0$ such that the following holds for all integers $n\geq n_0$. Let $\mathcal{G}=\{G_1,\ldots,G_n\}$ be a graph system of order $n$ and  {$\delta(\mathcal{G})\geq \frac{1-\alpha}{2}n$}. If $\mathcal{G}$ contains no transversal Hamilton cycles, then either $\mathcal{G}$ is $\kappa$-close to {$\mathcal{H}_{n-t}^t$} for some odd integer $t\in [n]$, or $\mathcal{G}$ is $\kappa$-close to a {half-split graph system}.
 \end{theorem}
 The proof of Theorem \ref{stability} was a combination of their newly developed transversal regularity-blow-up method \cite{2023chengBlowup} (see Theorem \ref{blowuplemma} in the next section)  and the absorption method. 
 As a continuation of the transversal version of Dirac's theorem (Theorem \ref{Hamilton-cycle}) and the stability result for transversal Hamilton cycles (Theorem \ref{stability}), Cheng, Sun, Wang and Wei \cite{wanting} characterized all graph systems with order $n$ and minimum degree at least $n/2-1$, each of which does not contain  transversal Hamilton cycles. 
 
Li, Li and Li \cite{2022liluyi} studied the existence of rainbow Hamilton paths in a graph system $\mathcal{G}=\{G_i:i\in[n]\}$ under Ore-type conditions. However, the edges of a rainbow Hamilton path only come from $n-1$ graphs of $\mathcal{G}$. Thus, it is worth studying whether $n-1$ graphs can guarantee the existence of a rainbow Hamilton path in a graph system $\mathcal{G}=\{G_i:i\in[n-1]\}$. It is well-known that Dirac's theorem implies that every $n$-vertex graph with minimum degree at least $(n-1)/2$ contains a Hamilton path. Cheng, Sun, Wang and Wei \cite{wanting} proved a transversal version of this result. 
\begin{theorem}[\cite{wanting}]\label{Hamilton-path}
     Let $n\in \mathbb{N}$  and $n\geq 3$. Suppose $\mathcal{G} = \{ G_1, \ldots , G_{n-1}\}$ is  a graph system of order $n$. If $\delta ( \mathcal{G} ) \geq
    \frac{n-1}{2}$, then $\mathcal{G}$ contains a transversal Hamilton path.
 \end{theorem}

Similar to the case of Hamilton cycles, 
Cheng and Staden \cite{2024cheng} established a stability result for transversal Hamilton paths;  Cheng, Sun, Wang and Wei \cite{wanting} characterized all graph systems with order $n$ and minimum degree at least $(n-3)/2$, each of which does not contain transversal Hamilton paths.

In 1971, Bondy \cite{Bondy1971} posed a  meta-conjecture which states that almost all  nontrivial sufficient conditions on a graph which {ensure} that the graph is Hamiltonian also {imply} that the graph is pancyclic, or there might  be a simple family of exceptional graphs. Some concepts similar to pancyclic have also attracted extensive attention, such as, vertex-pancyclic and panconnected. Here, we only illustrate those concepts in graph systems. Let $\mathcal{G}=\{G_1,\ldots,G_n\}$ be a system of graphs on a common vertex set $V$ of order $n$. We say 
\begin{itemize}
\item $\mathcal{G}$ is \textit{transversal  pancyclic} if there is a rainbow cycle of length $\ell$ in $\mathcal{G}$ for each integer $\ell\in [3,n]$; 
    \item $\mathcal{G}$ is \textit{transversal  vertex-pancyclic} if each vertex of $V$ is contained in a rainbow cycle of length $\ell$ for all $\ell\in [3,n]$; 
\item $\mathcal{G}$ is  \textit{rainbow panconnected} if for any two vertices $u,v\in V$, $\mathcal{G}$ contains a  rainbow path of length 
$\ell$ joining $u$ and $v$ for all $\ell\in [d_{\mathcal{G}}(u,v)+1,n]$, where $d_{\mathcal{G}}(x,y)$ is the length of the shortest rainbow path (if it exists) connecting $u$ and $v$ in $\mathcal{G}$.
\end{itemize}
Li, Li and Li \cite{2022liluyi,2024Li} obtained Dirac-type conditions for transversal pancyclicity, transversal vertex-pancyclicity and  rainbow panconnectedness respectively. The analogue results in bipartite graph systems were obtained by Bradshaw \cite{2021Bradshaw} and  Hu, Li, Li and Xu \cite{2024HuLi}. 
 

In addition to the Dirac-type condition, there are many other sufficient conditions to
guarantee the Hamiltonicity of a graph, such as Ore's condition \cite{1960Ore} and P\'{o}sa's  condition \cite{Posacond}. Bradshaw \cite{2021Bradshaw} and Li, Li and Li \cite{2022liluyi} proposed that it is also interesting to consider those conditions in graph systems. We list one of them here. 

\begin{problem}
[\cite{2021Bradshaw,2022liluyi}]
    Let $\mathcal{G} = \left \{ G_{1}, \ldots , G_{n}\right \}$ be a graph system of order $n$. If $\sigma_2 (G_i)\geq n$ for all $i\in [n]$, does $\mathcal{G}$ contain a transversal Hamilton cycle?
\end{problem}

\subsection{$F$-{factors}}
  Given two graphs $F$ and $G$, an \textit{$F$-tiling} is a set of vertex-disjoint copies of $F$ in $G.$ A \textit{perfect $F$-tiling} (or an \textit{$F$-factor}) in  $G$ is an $F$-tiling covering all vertices of $G$. 
  Characterizing sufficient conditions for the {existence of an $F$-factor}  is one of the central areas of extremal graph theory. The most celebrated one is the  Hajnal--Szemer\'edi theorem.
    \begin{theorem}[Hajnal--Szemer\'edi \cite{Hajnal-Szemeredi}, Corr\'adi--Hajnal \cite{1963Corradi} for $k=3$]\label{H-Z}
      Every $n$-vertex graph $G$ with $ n\in k\mathbb{N}$ and
      $ \delta(G)\geq(1-\frac{1}{k})n$ has a $K_{k}$-factor. Moreover, the minimum degree condition is sharp.
  \end{theorem}

For any $r$-vertex graph $F$, K\"{u}hn and Osthus \cite{KuhnOsthus} determined the smallest value of $\delta$ such that there is an  $n_0 \in \mathbb{N}$ and
a constant $c=c(F)$ such that any graph with $n\geq n_0$ vertices and minimum degree at least $\delta n+c$ contains an $F$-factor if $r| n$. For more results, please see the survey of  K\"{u}hn and Osthus \cite{KuhnOsthussurvey}.


Now, we survey some results on minimum degree conditions for the {existence of a} transversal 
$F$-factor in graph systems. Joos and Kim \cite{2021jooskim} established an exact minimum degree condition to guarantee the existence of a transversal perfect matching (i.e., $K_2$-factor) in graph systems. 
\begin{theorem}[\cite{2021jooskim}]
     Let $n\in \mathbb{N}$ and $n\geq 2$ $even.$ Suppose $\mathcal{G} = \{ G_1, \ldots , G_{\frac{n}{2}}\}$ is a graph system of order $n$ and  $\delta(\mathcal{G})\geq \frac{n}{2}$. Then $\mathcal{G}$ contains a transversal perfect matching.
\end{theorem}

 Bradshaw \cite{2021Bradshaw} proved the analogous  result in bipartite graph systems, which is
 equivalent to a special case of a matching theorem for $r$-partite $r$-uniform hypergraphs obtained  by Aharoni, Georgakopoulos and Spr\"ussel \cite{AhPM}, for the case $r = 3$. 


  Cheng, Han, Wang and Wang   \cite{2023chengspan} gave an asymptotic  version of the transversal Hajnal--Szemer\'{e}di theorem. 
  Furthermore, they proposed a conjecture for the  exact version of transversal Hajnal--Szemer\'{e}di theorem. 
\begin{conjecture}[\cite{2023chengspan}]
Let $\mathcal{G}=\{G_1,\ldots,G_{\frac{n}{k}\binom{k}{2}}\}$ be a graph system of order $n$. If $\delta(\mathcal{G})\geq (1-\frac{1}{k})n$, then $\mathcal{G}$ contains a transversal $K_k$-factor.
\end{conjecture}

Montgomery, M\"{u}yesser and Pehova \cite{2022Montgomery} showed that, in most cases, the same minimum degree bound is asymptotically sufficient for
the existence of transversal $F$-factors, which generalized the K\"{u}hn--Osthus theorem \cite{KuhnOsthus} on factors. For each graph $F$, let $\delta_F^f$ be the smallest real number $\delta\geq 0$ such that, for each
$\epsilon>0$ there is some $n_0$ such that, for every $n\geq n_0$ with $|V(F)|$ dividing $n$, if an $n$-vertex graph $H$ has minimum
degree at least $(\delta+\epsilon)n$, then $ H$ contains an $F$-factor. 
  \begin{theorem}[\cite{2022Montgomery}]\label{factor}
      Let $\varepsilon>0$ and let $F$ be a graph on $r$ vertices with $t$ edges.  If 
$\delta_{F}^f\geq\frac{1}{2}$
 or $F$ has a bridge, then let $\delta_{F}^{T}= \delta _{F}^f$, and otherwise let $\delta_{F}^{T}= \frac{1}{2}.$ Then, there is some $n_{0}$ such that the following holds for all $n\geq n_0.$ 
      
      Assume that $\mathcal{G}=\{G_1,\ldots,G_{tn}\}$ is a graph system of order $rn$ and $\delta(\mathcal{G})\geq(\delta_F^T+\varepsilon)rn$. Then, $\mathcal{G}$ contains a transversal $F$-factor.
  \end{theorem}
  Suppose $V_1,V_2,\ldots,V_k$ are disjoint vertex sets each of order $n$, and
$G$ is a $k$-partite graph on vertex classes $V_1,V_2,\ldots,V_k$. Denote the \textit{partite minimum degree} of $G$ by
$$\delta^{\prime}(G):=\min_{i\in[k]}\min_{v\in V_i}\min_{j\in[k]\setminus\{i\}}|N(v)\cap V_j|.
$$ 
Fischer \cite{1999Fischer} conjectured that the $k$-partite graph $G$ contains a $K_k$-factor if $\delta^\prime(G)\geq(1-1/k)n+1$ (actually his original conjecture did not include the term $+1$, but this stronger conjecture is known to be false). An approximate version of this conjecture was demonstrated independently by Keevash and Mycroft \cite{2015Keevash}, and by Lo and Markstr\"om \cite{2013Lo}. Later, Keevash and Mycroft \cite{Keevash2015} showed that the conjecture holds for any sufficiently large graph (in fact, they obtained a more general result for $r$-partite graphs with $r\geq k$). 
Cheng, Han, Wang and Wang  \cite{2023chengspan} extended the approximate result of Fischer's conjecture to the transversal version.
\begin{theorem}[\cite{2023chengspan}]\label{rainbow-HS}
    For every $\epsilon>0$ and integer $k$, there exists $n_0\in \mathbb{N}$ such that the following holds for all integers $n\geq n_0$. If $\mathcal{G}= \{ G_1, \ldots , G_{n\binom k2}\}$ is a system of $k$-partite graphs with a common vertex partition $ V_1, \ldots , V_k$\ each of size $n$ such that $\delta ^{\prime }( G_i) \geq ( 1- \frac 1k+ \varepsilon ) n$,  then $\mathcal{G}$ contains a transversal $K_{k}$-factor.
\end{theorem}

\subsection{Graphs with bounded bandwidth}
Montgomery, M\"{u}yesser and Pehova \cite{2022Montgomery} determined an asymptotically tight bound on
$\delta(\mathcal{G})$ for the existence of a transversal spanning tree with maximum degree at most $cn/ \log n$, generalizing the  tree-embedding theorem of Koml\'{o}s, S\'{a}rk\"{o}zy and Szemer\'{e}di \cite{2001komlosTree}.

It is noteworthy that many graphs, including $F$-factors, powers of Hamilton cycles, and trees, exhibit a somewhat bounded maximum degree and low connectivity. This property of low connectivity can be captured by the following notion of bandwidth. A graph $H$ has a \textit{bandwidth} at most $b$ if there exists an ordering $x_1,\cdots,x_n$ of $V(H)$ such that all edges $x_ix_j\in E(H)$ satisfy $|i-j|\leq b.$ B\"ottcher, Schacht and Taraz \cite{2009BottcherBW} proved the celebrated bandwidth theorem, which determines the asymptotically optimal minimum degree condition for a graph to contain a spanning subgraph $H$ with bounded maximum degree and low bandwidth.

Since the bandwidth theorem provides a common generalization of many classical results, it is natural to {give} such a generalization for the above transversal version results. Chakraborti, Im, Kim and Liu \cite{2023chakrabortiBandwidth} proved that
such a generalization holds by showing the following bandwidth theorem for graph transversals. 
\begin{theorem}[\cite{2023chakrabortiBandwidth}]\label{bandwidth}
    For every $\epsilon > 0$ and positive integers $\Delta, k,$ there exist $ \alpha > 0$ and $h_{0}> 0$ satisfying the following for every $h\geq h_0$. Let $H$ be an $n$-vertex graph with $h$ edges and bandwidth at most $\alpha n$ such that $\Delta (H) \leq \Delta$ and $\chi ( H) \leq k.$ If $\mathcal{G} = \{G_1, \ldots , G_h\}$ is a graph system of order $n$ and  $\delta(\mathcal{G})\geq\left(1-\frac1k+\varepsilon\right)n$, then $\mathcal{G}$ contains a transversal copy of $H$.
\end{theorem}

This covers several of the above results, and is asymptotically best possible for 
{many graphs $H$.} Note that by taking each $G_i$ in Theorem \ref{bandwidth}  to be the same graph, the bandwidth theorem of B\"{o}ttcher, Schacht and Taraz \cite{2009BottcherBW} can be deduced immediately. Therefore, Theorem \ref{bandwidth} generalizes and strengthens the original bandwidth theorem.


\section{Transversal via regularity}
The classical ``regularity-blow-up method'' has been used to prove many results concerning the embedding of a given spanning subgraph $H$ in a large graph $G$.
Such proofs are typically developed using two important tools: Szemer\'{e}di's regularity lemma \cite{regularity} and blow-up lemma of Koml\'{o}s, S\'{a}rk\"{o}zy and Szemer\'{e}di \cite{blowup}.  Cheng and Staden \cite{2023chengBlowup}  provided a tool needed to obtain transversal embeddings by applying the regularity-blow-up method. More precisely, they  proved a transversal blow-up lemma that embeds separable graphs into graph
systems. 

    
        
\begin{definition}[Regularity and superregularity] 
Suppose that $\mathcal{G}=\{G_c: c\in \mathscr{C}\}$ is a  system of bipartite graphs on a common vertex partition $V_1\cup V_2$ with color set $\mathscr{C}$. 
We say that 
\begin{enumerate}
    \item[\rm{\textbullet}] $\mathcal{G}$ is \textit{$(\varepsilon,d)$-regular} if whenever $V_i^{\prime}\subseteq V_i$ with $|V_i^{\prime}|\geq\varepsilon|V_i|$ for $i=1,2$ and $\mathscr{C}^\prime\subseteq\mathscr{C}$ with
$|\mathscr{C}^{\prime}|\geq\varepsilon|\mathscr{C}|$, we have

$$\left|\frac{\sum_{c\in\mathscr{C}'}|E(G_c[V_1',V_2'])|}{|\mathscr{C}'||V_1'||V_2'|}-\frac{\sum_{c\in\mathscr{C}}|E(G_c)|}{|\mathscr{C}||V_1||V_2|}\right|<\varepsilon $$
and $\sum_{c\in \mathscr{C}}|E(G_c)|
\geq d|\mathscr{C}||V_{1}||V_{2}|$.  
\item[\rm{\textbullet}] $\mathcal{G}$ is \textit{$(\varepsilon,d)$-superregular} if it is $(\varepsilon,d)$-regular and $\sum_{c\in\mathscr{C}}d_{G_c}(x)\geq d|\mathscr{C}||V_{3-i}|$ for all $x\in V_i$ with $i\in[2]$, and $|E(G_c)|\geq d|V_1||V_2|$ for all $c\in\mathscr{C}.$
\end{enumerate}

Note that if every $G_c$ with $c\in\mathscr{C}$ is the same, then $\mathcal{G}$ is $(\varepsilon,d)$-regular if and only if $G_c$ is $(\varepsilon,d)$-regular; and $\mathcal{G}$ is $(\varepsilon,d)$-superregular if and only if $G_c$ is $(\varepsilon,d)$-superregular.
\end{definition}
\begin{theorem}[Regularity lemma for graph systems, \cite{2023chengBlowup}]\label{regular}
For all integers $ L_{0}\geq 1$ and every $\varepsilon , \delta > 0$, there is an $n_{0}= n_{0}(\varepsilon , \delta , L_{0})$ such that for every $d\in [ 0, 1)$ and every graph system $\mathcal{G}= \{ G_{c}: c\in \mathscr{C} \}$ on vertex set $V$ of size $n \geq n_{0}$ with $\delta n\leq | \mathscr{C} | \leq \frac{n}{\delta} $,   there exists a partition of $V$ into $V_{0}, V_{1}, \ldots , V_{L}$, of $\mathscr{C}$ into $\mathscr{C}_0,\mathscr{C}_1,\ldots,\mathscr{C}_M$ and a spanning subgraph $G_c'$ of $G_c$ for each $c\in \mathscr{C}$  such that the following  properties hold:
\begin{enumerate}

\item [{\rm (i)}] $L_{0}\leq L, M\leq n_{0}$ and $| V_{0}| + | \mathscr{C}_{0}| \leq \varepsilon n;$
\item [{\rm (ii)}] $| V_{1}| = \cdots = | V_{L}| = | \mathscr{C}_{1}| = \cdots = | \mathscr{C}_{M}|;$
\item [{\rm (iii)}] $\sum _{c\in \mathscr{C} }d_{G_{c}^{\prime }}(v) > \sum _{c\in \mathscr{C} }d_{G_{c}}(v) - ( \frac{3d}{\delta ^{2}} + \varepsilon ) n^{2}$ for all $v\in V$ and $|E( G_{c}^{\prime})| > |E(G_{c})| - ( \frac{3d}{\delta ^{2}}+ \varepsilon ) n^{2}$ for all $c\in \mathscr{C};$
\item [{\rm (iv)}] if, for $c\in \mathscr{C}$, the graph $G_{c}^{\prime }$ has $an$ edge with both vertices in a single vertex cluster $V_{i}$ for some $i\in [ L]$,
then $c\in \mathscr{C}_0;$
\item [{\rm (v)}] for all triples $( \{ h, i\} , j)\in \binom {[L] }2 \times [M]$, we have that either $G_{c}^{\prime }[ V_{h}, V_{i}]$ = $\emptyset$ for all $c\in \mathscr{C} _{j}$, or
$\mathcal{G}_{hi, j}^{\prime }: = \{ G_{c}^{\prime }[ V_{h}, V_{i}] : c\in \mathscr{C}_{j}\}$ is $( \varepsilon , d)$-regular.
\end{enumerate}
The sets $V_i$ are called vertex clusters and the sets $\mathscr{C}_j$ are called color clusters,  while $V_0$ and $\mathscr{C}_0$ are the exceptional vertices and colors sets respectively.
\end{theorem}

\begin{definition}[Reduced graph system] 
Given a graph system $\mathcal{G}=\{G_c:c\in\mathscr{C}\}$ on $V$ and parameters $\varepsilon>0,\,d\in [0,1 )$ and $L_0 \geq 1$, the \textit{reduced graph system} $\mathcal{R}= \mathcal{R}( \varepsilon , d, L_0)$ 
is defined as follows. Apply Theorem \ref{regular} to $\mathcal{G}$ with parameters $\varepsilon,d,L_0$ to obtain a graph system $\mathcal{G}'=\{G'_c:c\in \mathscr{C}\}$, a partition $V_0,\ldots,V_L$ 
 of $V$ and $\mathscr{C}_0,\ldots,\mathscr{C}_M$ of $\mathscr{C}$,  where $V_0,\mathscr{C}_0$ are the exceptional sets, $V_1,\ldots,V_L$ are the vertex clusters and $\mathscr{C}_1,\ldots,\mathscr{C}_M$ are the color clusters. Then $\mathcal{R}=\{R_1,\ldots,R_M\}$ is a graph system consisting of $M$ graphs each on the same vertex set $[L]$, where $hi \in R_j$ if and only if  
 $\mathcal{G}_{hi,j}'$ is $( \varepsilon,d)$-regular for each $(\{ h, i\} , j) \in \binom {[ L] 
 }{2} \times [M]$.
\end{definition}

An $n$-vertex graph $H$ is \textit{$\mu$-separable} if there is $X\subseteq V(H)$ of size at most $\mu n$ such that $H-X$ consists of components of size at most $\mu n.$ With suitable small $\mu$,  separable graphs include many natural graphs, for example, $F$-factors for a fixed $F$, $2$-regular graphs, trees, powers of Hamilton cycles and graphs of small bandwidth.

\begin{theorem}[Transversal blow-up lemma, \cite{2023chengBlowup}]\label{blowuplemma}
      For all $\nu,d,\Delta,r>0$, where $r\geq 2$ is an integer, there exist $\epsilon,\mu,\alpha>0$ and $m_0\in \mathbb{N}$ such that the following holds for all integers $m\geq m_0$.  
    Assume $\mathcal{G}=\{G_i:i\in \mathscr{C}\}$ is a graph system with the following properties.

    \begin{itemize}
        \item There is a graph $R$ with vertex set $[r]$ and a partition $\mathscr{C}=\cup_{e\in E(R)}\mathscr{C}_e$ where $|\mathscr{C}_e|\geq \delta m$ for all $e\in E(R)$;
        \item for all $ij\in E(R)$ and $c\in \mathscr{C}_{ij}$, $G_c$ is bipartite with parts $V_i,V_j$, and $V$ is a vertex set of size $n$
with partition $V=V_1\cup V_2\cup \ldots \cup V_r$, where $m\leq |V_i|\leq \frac{m}{\delta}$ for each $i\in [r]$;
\item for all $ij\in E(R)$,
\begin{itemize}
    \item for all ${V_h}'\subseteq V_h$ $(h=i,j)$ and $\mathscr{C}_{ij}'\subseteq \mathscr{C}_{ij}$  with $|{V_h}'|\geq \epsilon |V_h|$ $(h=i,j)$ and $|\mathscr{C}_{ij}'|\geq \epsilon |\mathscr{C}_{ij}|$, we have that
$$
\sum_{c\in \mathscr{C}_{ij}'}|E(G_c[{V_i}', {V_j}'])|\geq d|\mathscr{C}_{ij}'||{V_i}'||{V_j}'|;
$$
\item for every $v\in V_h$ with $h=i,j$ we have $\sum_{c\in \mathscr{C}_{ij}'}d_{G_c}(v)\geq d|\mathscr{C}_{ij}'|m$, and for every $c\in \mathscr{C}_{ij}$, we have $|E(G_c)|\geq dm^2$.
\end{itemize}
\end{itemize}

Suppose that $H$ is a graph with the following properties. 
\begin{itemize}
    \item $\Delta(H)\leq \Delta$;
    \item $H$ is $\mu$-separable;
    \item $H$ has vertex partition $A_1\cup A_2\cup \ldots \cup A_r$ such that $|A_i|=|V_i|$ for all $i\in [r]$, and for every $xy\in E(H)$
there is $ij\in E(R)$ such that $x\in A_i$ and $y\in A_j$;
\item $|E(H[A_i,A_j])|=|\mathscr{C}_{ij}|$  for all $ij\in E(R)$;
\item for $i \in [r]$, there is a set  $U_i\subseteq A_i$ with $|U_i|\leq \alpha |A_i|$ and for each $x\in U_i$, there exists a target set $T_x\subseteq V_i$ with $|T_x|\geq \nu m$.
\end{itemize}
   Then $\mathcal{G}$ contains a transversal copy of $H$ such that every $x\in U_{i}$ is embedded inside $T_{x}$ for $i\in[r]$.
 \end{theorem}
By utilizing the transversal blow-up lemma, Cheng and Staden  \cite{2023chengBlowup} proved a stability result for transversal Hamilton cycles,  and characterized a large class of spanning $3$-uniform
linear hypergraphs $H$ such that any sufficiently large uniformly dense $n$-vertex $3$-uniform hypergraph with minimum vertex degree $\Omega(n^2)$ contains $H$ as a subhypergraph.

\section{Random graph systems}
In this section, we survey some results on the existence of { transversal structures} 
in a system of random graphs in various settings. 
Let $0\leq p\leq1$ and $G(n,p)$ be the binomial \textit{random graph}. For graphs $G=(V,E)$ and $G'=(V',E')$, we set $G\cap G'= (V\cap V', E\cap E')$. 
Given an $n$-vertex graph $G$, we call $G(p):=G\cap G(n,p)$ the \textit{random subgraph} of $G$. In particular, $G(p)=G(n,p)$ when $G$ is the complete graph $K_n.$ 
We say that $G(p)$ has a graph property $\mathcal{P}$ with high probability, or whp for brevity, if the probability that $G(p)$ has $\mathcal{P}$ approaches $1$ as $n$ tends to infinity. 
Sudakov \cite{SudakovRobust} pointed out that 
{$G(p)$} can be used to measure the robustness of graph properties.


A natural question concerns the existence of transversal spanning structures in random graph systems. 
A graph property is called \textit{monotone increasing} if it is closed under the addition of
edges. 
For a monotone increasing graph property, the \textit{resilience}  quantifies the robustness in terms of the number of edges one must delete from a graph, locally or globally, in order to destroy the property. 
Ferber, Han and Mao \cite{2022FerberRandom}  proved a resilient result for a system of $n$ graphs with order $n$ where each graph is an independent copy of $G(n, p)$ provided that $np = \omega(\log n)$, {and they also noted that the bound for $p$ is probably not tight.}

\begin{theorem}[\cite{2022FerberRandom}]\label{ferberRD}
    For every $\epsilon>0$, there exists an absolute constant $c= c(\epsilon)$ such that for any sufficiently large $n$ and $p\geq \frac{c\log n}n$ the following holds. Let $\mathcal{G}=\{G_1,\ldots,G_n\}$ be a system of $n$ independent copies of $G(n,p)$ on the same vertex set $V$ of size $n$. For all $i\in [n]$, let $H_i\subseteq G_i$ be such that $d_{H_i}(v)\geq (\frac{1}{2}+ \epsilon)d_{G_i}(v)$ for all $v\in V$. Then, whp $\{ H_{1}, \ldots , H_{n}\}$ contains a  transversal Hamilton cycle.
\end{theorem}

{A similar result for transversal perfect matchings} was also obtained in \cite{2022FerberRandom}. Note that the above theorem implies that if $\mathcal{G}=\{G_1,\ldots,G_n\}$ is a graph system of order $n$ with $\delta(\mathcal{G})\geq (1/2+ \epsilon)n$, and $\mathcal{F}=\{F_1,\ldots,F_n\}$ where $F_i=G_i(p)$  for some $p\geq {c\log n}/n$, then whp $\mathcal{F}$ contains a  transversal Hamilton cycle. 
  Anastos and Chakraborti \cite{2023Anastosrobust} improved this  by removing the error term in the minimum degree condition with spread techniques.
\begin{theorem}[\cite{2023Anastosrobust}]\label{Anastos}
    {There exists an absolute $c$ such that the following holds.
Suppose that $\mathcal{G} =\{G_1,\ldots,G_n\}$ is a graph system of order  $n$ with $\delta(\mathcal{G})\geq \frac{n}{2}$, and $p \geq \frac{c\log n}n$. Let $\mathcal{F}=\{F_1,\ldots,F_n\}$ where $F_i=G_i(p)$. Then, whp $\mathcal{F}$  contains a transversal Hamilton cycle.}
\end{theorem}


In the above theorem, each graph in the given system intersects with the same random graph $G(n,p)$. They \cite{2023Anastosrobust} also considered a system that consists of $n$ independent random subgraphs of a Dirac graph.

\begin{theorem}[\cite{2023Anastosrobust}]\label{rthc}
    There exists an absolute constant $c$ such that the following holds. Suppose $G$ is 
    a graph on $n$ vertices
    with $\delta(G)\geq \frac{n}{2}$, and $p\geq \frac{c\log n}{n^2}.$ Let  $\mathcal{G}=\{G_1,\ldots,G_n\}$ be a family of graphs where each $G_i$ is independently distributed as $G(p)$, for $i\in [n]$. Then, whp  $\mathcal{G}$ contains a transversal Hamilton cycle.
\end{theorem}
{In Theorem \ref{Anastos} and Theorem \ref{rthc}, 
the lower bound on $p$ is optimal up to a constant factor $c$.} Naturally, they \cite{2023Anastosrobust} hoped for a combined generalization of the above two theorems.
\begin{problem}[\cite{2023Anastosrobust}]\label{2023Anastos}
    Does there exist a constant $c$ such that the following holds? Let $\mathcal{G}=\{G_1,\ldots,G_n\}$ 
    be a graph system of order $n$
    with $\delta(\mathcal{G})\geq \frac{n}{2}$, and $p\geq \frac{c\log n}{n^2}$. For each $i\in [n]$, let $F_i=G_i(p)$ and let $\mathcal{F}=\{F_1,\ldots,F_n\}$. Then, whp $\mathcal{F}$ contains a transversal Hamilton cycle.
\end{problem}
{ In the same paper, they \cite{2023Anastosrobust} showed that Problem \ref{2023Anastos} has a positive answer when $n$ is odd}, whereas it  turns out to be false when $n$ is even.

Han, Hu and Yang \cite{2023HanH-S} 
obtained a transversal robust version of the Hajnal--Szemer\'{e}di theorem \cite{Hajnal-Szemeredi}, {where the bound for $p$ is not tight. }
\begin{theorem}[\cite{2023HanH-S}]
    {Let $r\geq 3$ be an integer and let $\gamma > 0.$ There exists a constant $c= c( r, \gamma )$ such that for any sufficiently large $n\in r\mathbb{N}$ and $p\geq cn^{- \frac{2}{r}}( \log n) ^{1/\binom{r}{2}}$ the following holds. Suppose that $\mathcal{G} =\{G_1,\ldots,G_{\frac nr\binom r2}\}$ is a graph system of order $n$ with $\delta(\mathcal{G})\geq ( 1- \frac{1}{r}+ \gamma ) n$. Let $\mathcal{F}=\{F_1,\ldots,F_{\frac nr\binom r2}\}$ where $F_i=G_i(p)$.  Then whp $\mathcal{F}$ contains a transversal $K_{r}$-factor.}
\end{theorem}


The transversal robust version of the multipartite Hajnal--Szemer\'{e}di theorem \cite{Hajnal-Szemeredi} was also obtained in \cite{2023HanH-S}.

\section{Digraph systems}
In this section, we consider transversal problems in digraph systems. Let $D$ be a \textit{digraph}  with no loops and at most one edge in each direction between every pair of vertices. 
An \textit{oriented graph} is a digraph with no cycle of length two. 
Let $\delta^{+}(D)$ and $\delta^{-}(D)$ be the \textit{minimum out-degree} and \textit{minimum in-degree}, respectively. Let $\delta^0(D):=\min\{\delta^+(D),{\delta^-}(D)\}$ be the \textit{minimum semi-degree} of $D.$  Let $\mathcal{D} =\{D_1,D_2,\dots,D_m\}$ be a system of not necessarily distinct digraphs with common vertex set $V$. Define $\delta^+(\mathcal{D}):=\min\{\delta^+(D_i):i\in [m]\}$ and $\delta^0(\mathcal{D}):=\min\{\delta^0(D_i):i\in [m]\}$. 

A \textit{tournament} is an oriented graph where every pair of distinct vertices are adjacent. 
We call a tournament $T$ \textit{strongly connected} if for every pair of vertices $x,y$ in $T$, there is a directed path from $x$ to $y$. 
We say that a digraph $D$ is \textit{transitive} if $xy$ and $yz$ are arcs in $D$ with $x\neq z$, then the arc $xz$ is also in $D$.

 Let $T_k$ be the transitive tournament on $k$ vertices and $\mathcal{T}_k$ be the family of all tournaments on $k$ vertices.
 Czygrinow, DeBiasio, Kierstead and Molla \cite{Czygrinow} proved that every $n$-vertex  digraph $D$ with $\delta^+(D)\geq (1-1/k)n$ contains
a $T_k$-factor for integers $n,k$ with $k|n$. Treglown \cite{Treglown} proved that for  
$T\in \mathcal{T}_k$ with $k\geq 3$ and sufficiently large integer $n$ with $k|n$, every digraph $D$ on $n$ vertices with $\delta^0(D)\geq (1-1/k)n$ contains
a $T$-factor.
 
 {In light of} Theorem \ref{rainbow-HS}, it is natural to seek an  analogue of the transversal Hajnal--Szemer\'{e}di theorem \cite{Hajnal-Szemeredi} in { digraph systems}. Cheng, Han, Wang and Wang \cite{2023chengspan} established the minimum out-degree and semi-degree condition for the existence of  transversal tournament factors in digraph systems. 
\begin{theorem}[\cite{2023chengspan}]
    For every integer $ k\geq 3$, $T\in \mathcal{T}_k$ and real $\varepsilon > 0$,  there exists $n_{0}\in \mathbb{N}$ such that the following holds for all integers $n\geq n_0$ and $n\in k\mathbb{N}$. If $\mathcal{D}=\left\{D_1,\ldots,D_m\right\}$, $m=\frac nk\binom k2$, is a system of $n$-vertex digraphs on the same vertex set such that $\delta ^{+}(\mathcal{D}) \geq ( 1- \frac 1k+ \varepsilon ) n$ (resp. $\delta^0(\mathcal{D}) \geq (1-\frac{1}{k}+\varepsilon)n$), then $\mathcal{D}$ contains a transversal  $T_k$-factor 
    (resp. $T$-factor).
\end{theorem}

It is well-known that every tournament contains a Hamilton path, and every strongly connected tournament contains a Hamilton cycle.  Chakraborti, Kim, Lee and Seo \cite{2023Tournament} proved the following result concerning the existence of transversal directed Hamilton paths and cycles in a system of tournaments. 

\begin{theorem}[\cite{2023Tournament}]
     For every sufficiently large $n$, 
     \begin{enumerate}
        
     \item [\rm (i)] every system of $(n-1)$ tournaments with order $n$ contains a transversal directed Hamilton path;
     \item [\rm(ii)] every system $\mathcal{T}$ of $n$ tournaments with  order $n$ satisfies the following. If all tournaments in $\mathcal{T}$ possibly except one are strongly connected, then $\mathcal{T}$ contains a transversal directed Hamilton cycle.
     \end{enumerate}
\end{theorem}

Thomason \cite{1986ThomasonTour} proved that every sufficiently large tournament contains Hamilton paths and cycles with all possible orientations, except possibly the consistently oriented Hamilton cycle. Chakraborti, Kim, Lee and Seo \cite{2024KimTour} established the transversal generalization of this result.
\begin{theorem}[\cite{2024KimTour}]
 For every sufficiently large $n$, 
 \begin{enumerate}
   
 \item [\rm (i)]every system  of $(n-1)$ tournaments with order $n$ contains a transversal of every orientation of the Hamilton path;
 
 \item [\rm (ii)] every system of $n$ tournaments with order $n$ contains a transversal of every orientation of the Hamilton cycle except possibly the directed one.
    \end{enumerate}
\end{theorem}

Babi{\'n}ski, Grzesik and Prorok \cite{2023BabinskiTria} studied the transversal version of Mantel's theorem in digraph systems  and oriented graph systems. More precisely, they determined the minimum number of edges in each graph that {ensures} the existence of transversal transitive triangles or transversal directed triangles in digraph systems and oriented graph systems. Gerbner, Grzesik, Palmer and  Prorok \cite{2024GerbnerDstar} 
determined the minimum number of edges {required in each digraph of a digraph system that guarantees the existence of any rainbow directed star.}

Ghouila-Houri \cite{Houri} showed that any $n$-vertex digraph with  minimum semi-degree at least $n/2$ contains a directed Hamilton cycle. Cheng, Li, Sun and Wang \cite{Chenglisun} provided a transversal generalization of Ghouila-Houri's theorem, which solved a problem proposed by Chakraborti, Kim, Lee and Seo \cite{2024KimTour}. 

\begin{theorem}[\cite{Chenglisun}]
   Let $\mathcal{D} = \left \{ D_{1}, \ldots , D_{n}\right \}$ be a system of digraphs on a  common vertex set $V$ of sufficiently large order $n.$ If $\delta^0 (\mathcal{D})\geq \frac{n}{2}$, then $\mathcal{D}$ contains a transversal directed  Hamilton cycle.
\end{theorem}

DeBiasio and Molla \cite{ADHC} proved the anti-directed Hamilton cycle (an orientation in which consecutive edges  alternate  direction) is guaranteed to appear if $\delta^0(D) \geq n/2+1$. 
DeBiasio, K\"uhn, Molla, Osthus and Taylor \cite{AOHC} showed that for sufficiently large $n$, every $n$-vertex digraph $D$ with $\delta^0 (D)\geq n/2$ contains every orientation of a Hamilton cycle except, possibly, the anti-directed one. It {would be}  interesting to generalize the above two results to the  {transversal setting}.

Keevash,  K\"uhn and Osthus \cite{2009KeevashOriented} proved the oriented version of Dirac's theorem, which states that { for sufficiently large $n$}, every $n$-vertex digraph with minimum semi-degree at least $(3n-4)/8$ contains a directed Hamilton cycle. Considering a transversal version  of this theorem would be interesting. 

\begin{problem}[\cite{Chenglisun}]
    Let $\mathcal{D} = \left \{ D_{1}, \ldots , D_{n}\right \}$ be a system of oriented graphs with common vertex set $V$ of size $n.$ {For sufficiently  large $n$,} if $\delta^0 (\mathcal{D})\geq \lceil\frac{3n-4}{8}\rceil$, does $\mathcal{D}$ contain a transversal directed  Hamilton cycle?
\end{problem}

\section{Hypergraph systems}

In this section, we will survey results in $k$-uniform hypergraph systems. For a set $V$, we denote by $\binom{V}{k}$ the set of subsets of size $k$ of $V$. A \textit{$k$-uniform hypergraph} (or \textit{$k$-graph}) ${H}=(V,E)$ consists of a vertex set $V$ and an edge set $E$ which is a family of $k$-element subsets of $V$, i.e., $E\subseteq\binom Vk.$ For any $S\subseteq V$, the degree of $S$ in ${H}$, denoted by $\deg_{H}(S)$, is the number of edges containing $S.$  For any integer $s\geq 0$, define the \textit{minimum $s$-degree} $\delta_s({H})$ to be $\min\{\deg_{H}(S):S\in \binom{V({H})}{s}\}$. Let $\mathcal{H} =\{H_1,H_2,\dots,H_m\}$ be a system of not necessarily distinct $k$-graphs with common vertex set $V$. Define $\delta_s(\mathcal{H})=\min\{\delta_s(H_i):i\in [m]\}$. 

The problem of characterizing $k$-graphs which contain a Hamilton cycle or a perfect matching is perhaps the most prominent one in hypergraph theory. In the remainder of this section, we focus on those two most well studied structures.


\subsection{Hamilton cycle}

We say that a $k$-graph ${C}$ is an \textit{$\ell$-cycle} if there exists a cyclic ordering of the vertices of ${C}$ such that every edge of ${C}$ consists of $k$ consecutive vertices and every two consecutive edges (in the natural ordering of the edges) intersects in precisely $\ell$ vertices. We say that $C$ is a \textit{Hamilton $\ell$-cycle} of a  $k$-graph ${H}$ if $C$ is a spanning sub-$k$-graph of $H$. We call $C$ a  \textit{tight Hamilton cycle} if $\ell =k-1$, a \textit{loose Hamilton cycle} if $\ell =1$.  A Hamilton cycle of a $k$-graph $H$ is a cycle which is a spanning subhypergraph of $H$. Note that if a $k$-graph ${H}$ on $n$ vertices contains a Hamilton $\ell$-cycle then $(k-\ell)|n$, since every edge of the cycle contains exactly $k-\ell$ vertices which were not contained in the previous edge.   

 In 1999, Katona and Kierstead \cite{1999KatonaHCHyper} proved a Dirac-type theorem in hypergraphs: every $n$-vertex $k$-graph $H$  with $\delta_{k- 1}( {H} ) > ( 1- \frac 1{2k}) n+ 4- k- \frac 5{2k}$ contains a tight Hamilton cycle. 
They anticipated that a much stronger result is true.
\begin{conjecture}[\cite{1999KatonaHCHyper}]\label{conjtightHC}
    Let $H$ be a $k$-graph on $n\geq k+1\geq 4$ vertices. If $\delta_{k-1}(H)\geq \lfloor\frac{n-k+3}{2}\rfloor$, then $H$ has a tight Hamilton cycle.
\end{conjecture}
R\"odl, Ruci\'nski and Szemer\'edi  \cite{2006RodlHCHyper,2008RodlHCHyper} confirmed this asymptotically, and they \cite{2011RodlHCHyper} also gave the exact version for $k = 3$.
\begin{theorem}[\cite{2006RodlHCHyper,2008RodlHCHyper,2011RodlHCHyper}]\label{2008Rodl}
    Let $k\geq 3$, $\gamma > 0$ and ${H}$ be an $n$-vertex $k$-graph, where $n$ is sufficiently large.  If $\delta_{k-1}( {H}) \geq ( \frac{1}{2}+ \gamma ) n$, then ${H}$ contains a tight Hamilton cycle. Furthermore, when $k= 3$ it is  enough to have $\delta_{2}({H}) \geq \lfloor \frac{n}{2}\rfloor$.
\end{theorem}
However, K\"uhn and Osthus \cite{2014K-OHCHyper} and Zhao \cite{2016ZhaoHCHyper} noted that it is much more difficult to determine the minimum $d$-degree condition for tight Hamilton cycle for $d\in[k-2].$ 
The best general bound for such a problem was given by Lang and Sanhueza-Matamala \cite{2022Lang}, {and Polcyn}, Reiher, R\"odl and Sch\"ulke \cite{2021Polcyn} independently. 
\begin{theorem}[\cite{2022Lang,2021Polcyn}]\label{2022Lang}
    Let $k\geq 3$, $\gamma > 0$ and ${H}$ be an $n$-vertex $k$-graph,  where $n$ is sufficiently large. If  $\delta _{k-2}({H}) \geq ( \frac{5}{9}+ \gamma ) \binom n2$, then ${H}$ contains a tight Hamilton cycle. 
\end{theorem}

Let $1\leq d,\ell\leq k-1$. For $n\in (k-\ell)\mathbb{N}$, define $h_d^\ell(k,n)$ to be the smallest integer $h$ such that
every $n$-vertex $k$-graph $H$ satisfying $\delta_d(H)\geq h$ contains a Hamilton $\ell$-cycle. Han and Zhao \cite{hanjie} 
proved that
\begin{align}\label{eq:degree}
    h_d^{k-1}(k,n)\geq \left(1-\binom{t}{\lfloor\frac{t}{2}\rfloor}\frac{\lceil\frac{t}{2}\rceil^{\lceil\frac{t}{2}\rceil}(\lfloor\frac{t}{2}\rfloor+1)^{\lfloor\frac{t}{2}\rfloor}}{(t+1)^t+o(1)}\right)\binom{n}{t},
\end{align}
where $d\in [k-1]$ and $t=k-d$. In particular, $h_{k-3}^{k-1}(k,n)\geq (5/8+o(1))\binom{n}{3}$. Lang and Sanhueza-Matamala \cite{2022Lang} proved $h_d^{k-1}(k,n)\leq 2^{-1/(k-d)}$ and they conjectured that the minimum $d$-degree threshold for $k$-uniform tight Hamilton cycles coincides with the lower bounds in \eqref{eq:degree}.

Recently, Cheng, Han, Wang, Wang and Yang \cite{2021ChengHyp} and Tang, Wang, Wang and Yan \cite{2023TangHyp} extended  Theorem \ref{2008Rodl} and Theorem \ref{2022Lang} to transversal {versions}, respectively.

\begin{theorem}[\cite{2021ChengHyp,2023TangHyp}]\label{tightTHC}
    For every $k\geq3,\gamma>0$, there exists $n_{0}$ such that the following holds for $n\geq n_{0}.$ Let $\mathcal{H}=\{H_1,\ldots,H_n\}$ be a $k$-graph system of order $n$.
    \begin{enumerate}
    \item[{\rm (i)}] If $\delta_{k-1}(\mathcal{H})\geq(\frac{1}{2}+\gamma)n$, then $\mathcal{H}$ admits a transversal tight Hamilton cycle.
        \item[{\rm (ii)}] If $\delta_{k-2}(\mathcal{H})\geq(\frac{5}{9}+\gamma)\binom n2$, then $\mathcal{H}$ admits a transversal tight Hamilton cycle.
    \end{enumerate}
\end{theorem}


Cheng, Han, Wang, Wang and Yang  \cite{2021ChengHyp} suspected that the threshold for transversal tight Hamilton cycles in a $k$-graph system
is the same {as} the threshold conjectured for tight Hamilton cycles in a single $k$-graph (i.e., Conjecture \ref{conjtightHC}) and provided the following conjecture.

\begin{conjecture}[\cite{2021ChengHyp}]
     Suppose $\mathcal{H}=\{H_1,\ldots,H_n\}$ is a $k$-graph system of order $n\,(n\geq k+1 \geq 4)$ and $\delta_{k-1}(\mathcal{H}) \geq \lfloor \frac{n-k+3}{2}\rfloor$, then $\mathcal{H}$ contains  a transversal tight Hamilton cycle.
\end{conjecture}

 {Motivated} by \eqref{eq:degree}, Tang, Wang, Wang and Yan \cite{2023TangHyp} proposed the following conjecture, which concerns the minimum $d$-degree condition for transversal tight Hamilton cycles.
\begin{conjecture}[\cite{2023TangHyp}]
    For every $k\geq d+1$, $\mu >0$, there exists an integer $n_{0}$ such that the following holds for $n \geq n_{0}.$ Let $\mathcal{H}=\{H_1,\ldots,H_n\}$ be a $k$-graph system of order $n$.
    \begin{enumerate}
        \item[{\rm (i)}] If $k\geq 4$ and  $\delta_{k-3}(\mathcal{H})\geq(\frac{5}{8}+\mu)\binom n3$, then $\mathcal{H}$ admits a transversal tight Hamilton cycle.
        \item[{\rm (ii)}] If $\delta _d(\mathcal{H}) \geq h_d^{k-1}( k, n) + \mu \binom nd$, then $\mathcal{H}$ admits a transversal tight  Hamilton cycle.
    \end{enumerate}
\end{conjecture}


Furthermore, the problem of giving sufficient conditions for  transversal Hamilton $\ell$-cycles with $\ell \in [k-2]$ is  still open. It would be very nice to settle this problem. 


Theorem \ref{tightTHC} implies that if we consider $(k-1)$-degree and $(k-2)$-degree conditions, then tight Hamilton {cycles are} color-blind.  
Gupta, Hamann, M\"{u}yesser, Parczyk and Sgueglia \cite{2023GuptaPowerH} provided a unified approach to find color-blind graphs, 
which is general enough to recover many previous results and obtain novel transversal variants of several classical Dirac-type results. The precise statement of their main result is quite technical, we refer the reader to see \cite[Theorem 2.6]{2023GuptaPowerH}. The following theorem
lists some applications in their paper and they believed that their setting can capture even more families of hypergraphs.



\begin{theorem}[\cite{2023GuptaPowerH}]\label{Gupta}
    Under $d$-degree condition, the following families of hypergraphs are all color-blind.
    \begin{enumerate}
\item[{\rm (i)}] The family of the r-th powers of Hamilton cycles for $d=1$ and fixed $r\geqslant2$.
\item[{\rm (ii)}] The family of $k$-uniform Hamilton $\ell$-cycles for the following ranges of $k,\, \ell$ and $d$.
\begin{itemize}
\item $1<\ell<\frac{k}{2}$ and $d=k-2;$
\item $1\leq\ell<\frac{k}{2}$ or $\ell=k-1,and$ $d=k-1;$ 
\item $\ell=\frac{k}{2}$ and $\frac{k}{2}<d\leq k-1$ with $k$ even.
\end{itemize}
 \end{enumerate}
\end{theorem}


\subsection{Matching and factor}

A {\textit{$t$-matching}} $M$ in a hypergraph $H$ is a subset of $E(H)$ consisting of $t$ pairwise disjoint edges, which is \textit{perfect} if $V(M) = V(H)$. We denote by $\binom{[n]}{r}$ the set of subsets (also called edges) of size $r$ of $[n].$ For a set $X$, let $2^X$ be the set of all subsets of $X$. 
Let ${F}$ be a subfamily of $2^X$, denote by $\nu({F})$ its \textit{matching number}, that is, the largest $m$ such that ${F}$ contains $m$ pairwise disjoint sets. A classical problem in extremal set theory is to determine the maximum size of ${F}$ with  fixed $\nu({F})$. In 1965, Erd\H{o}s \cite{1965erdosmatch} proposed the following conjecture, which is known as Erd\H{o}s matching conjecture.

  \begin{conjecture}[Erd\H{o}s matching conjecture, \cite{1965erdosmatch}]\label{EMC}
    If $n\geq(s+1)k$  and a family ${F}\subset\binom{[n]}k$ satisfies
$\nu ( {F} ) \leq s$, then $| {F} | \leq \max\{\binom{n}{k}-\binom{n-s}{k},\binom{(s+1)k-1}k\}.$
\end{conjecture}
Notice that the  bound in Conjecture \ref{EMC} is tight. In fact, one may consider the edge set of the complete $k$-graph on $(s+1)k-1$ vertices or  the $k$-graph on $n$ vertices in which every edge intersects a fixed set of $s$ vertices. Setting $s=1$ in Conjecture \ref{EMC}, then it corresponds to the well-known  Erd\H{o}s--Ko--Rado theorem \cite{EKR}. 
There are numerous papers devoted to the Erd\H{o}s matching conjecture. For recent development on this conjecture, see \cite{2012AlonHypmatch,2013FranklHypmatch,2016HanHypmatch,2023KolupaevHypmatch}. 

Many scholars have attempted to generalize the Erd\H{o}s matching conjecture to {the} transversal version in hypergraph systems. Matsumoto and 
 Tokushige \cite{1989MatsumotoEKR} generalized the  Erd\H{o}s--Ko--Rado theorem  to two {families.}

\begin{theorem}[\cite{1989MatsumotoEKR}]\label{1989MatsumotoEKR}
Let 
$n,k,\ell$ be integers such that $n\geq \max\{2k,2\ell\}$. Assume that $\mathcal{H}_1\subseteq \binom{[n]}k$ and $\mathcal{H}_2\subseteq \binom{[n]}\ell$. If $|\mathcal{H}_1||\mathcal{H}_2|>\binom{n-1}{k-1}\binom{n-1}{\ell-1}$, then there exist two disjoint subsets  $X_1\in \mathcal{H}_1, X_2\in \mathcal{H}_2$.
\end{theorem}

It is natural to extend {the} Erd\H{o}s--Ko--Rado theorem to more than two hypergraphs. For this, Aharoni and Howard \cite{2011AhaRMconj} and Huang, Loh and Sudakov \cite{2011HuangMatch} proposed the following conjecture independently.

\begin{conjecture}[\cite{2011AhaRMconj,2011HuangMatch}]\label{2011} 
    Let $\mathcal{H} = \{ H_{1}, \ldots , H_{t}\}$ be a system of $k$-graphs on a common vertex set of size $n$. If
$|E(H_i)|>\max\left\{{\binom{n}{k}}-{\binom{n-t+1}{k}},{\binom{kt-1}{k}}\right\}$
for  all $1\leq i\leq t$,  then $\mathcal{H}$ has a transversal $t$-matching.
\end{conjecture}
When $k = 2$, this conjecture was confirmed by Meshulam (see \cite{2011AhaRMconj}).
Indeed, the case $k=2$ of this conjecture follows directly from an earlier result established by Akiyama and Frankl \cite{1985Akiyama}, as restated in \cite{1987Franklshift}.
Huang, Loh and Sudakov  \cite{2011HuangMatch} proved that Conjecture \ref{2011} holds for $n > 3k^2t$. 
Keller and Lifshitz \cite{2011kellerAhaconj} proved that Conjecture \ref{2011} holds when $n\geq f(t)k$ for some large constant $f(t)$ which only depends on $t$, and this was further improved the bound on $n$ to $12kt\log(e^2t)$ by Frankl and
Kupavskii \cite{2020FranklPM}. Lu, Wang and Yu \cite{2023LuMatch} improved this bound to $n >2kt$ for sufficiently large $t$. By using the sharp threshold techniques developed in \cite{2024Keevashfun}, Keevash, Lifshitz, Long and Minzer \cite{2019keevash} proved a more general result with $n=\Omega(kt)$, which is the optimal order of magnitude.  Kupavskii \cite{2021KupavskiiMatch} gave the concrete dependencies on the parameters by showing the conjecture holds for $n>3ekt$ with $t>10^7$. Recently, Gao, Lu, Ma and Yu \cite{2022GaoHypmatch} showed that Conjecture \ref{2011} is true when $k = 3$ and $n$ is sufficiently large. 

Aharoni and Howard \cite{2011AhaRMconj} proposed the following conjecture in $k$-partite $k$-graph systems.
\begin{conjecture}[\cite{2011AhaRMconj}]\label{conj-Aha-2}
    Let $\mathcal{H}=\{H_1,\ldots,H_t\}$ be a $k$-partite $k$-graph system with $n$ vertices in each partition class. If  
    $|E(H_i)|\geq (t-1)n^{k-1}$ for all $i\in [t]$, then $\mathcal{H}$ admits a transversal $t$-matching. 
\end{conjecture}
Aharoni and Howard  proved that Conjecture \ref{conj-Aha-2} holds when $t=2$ \cite{2011AhaRMconj} or $k\leq 3$ \cite{2017AhaEKR}. Lu and Yu \cite{2018LuHRP} proved Conjecture \ref{conj-Aha-2} for $n >3(k-1)(t-1)$. Kiselev and Kupavskii \cite{KiselevHRM} showed that this conjecture is true for $t\geq t_0$ where $t_0$ is a constant.

The degree conditions for the existence of large matchings in uniform hypergraphs have been extensively studied. 
K\"{u}hn, Osthus and Treglown \cite{2013Kuhn3-uniform} (also see Khan \cite{2013Khan3-uniform}) and Khan \cite{2016Khan4-uniform} determined the minimum $1$-degree condition for an $n$-vertex $3$-graph and $4$-graph { to contain} a perfect matching, respectively. Lu, Yu and Yuan \cite{2021LuHPM3-uniform} and  Lu, Wang and Yu \cite{2022LuHPM4-uniform}  proved the transversal version of the above results.
\begin{theorem}[\cite{2022LuHPM4-uniform,2021LuHPM3-uniform}]
    Let $n$ be a sufficiently large integer and $k\in \{3,4\}$ with $k|n$. If  $\mathcal{H}=\{H_1,\ldots,H_{\frac{n}{k}}\}$ is a $k$-graph system of order $n$  with $\delta_1(\mathcal{H})>\binom{n-1}k-\binom{(k-1)n/k}k$, then  $\mathcal{H}$ contains a transversal perfect matching. 
\end{theorem}



Recently, Lu and Wang \cite{2024Lu3-3match} determined the minimum $1$-degree condition that guarantees a $3$-partite $3$-graph system admits transversal perfect matchings. This generalized a result of Lo and Markstr\"om \cite{2014LoPM3-3} on the $1$-degree threshold for the existence of perfect matchings in $3$-partite $3$-graphs. Moreover, in their proof, they used fractional transversal matching theory obtained by Aharoni, Holzman and Jiang \cite{2019Aharonifrac}. 

It is well-known that perfect matchings are closely related to its fractional counterpart. Given a $k$-graph $H$, a \textit{fractional matching} is a function $f: E( H) \to [ 0, 1]$, subject to the requirement that $\sum_{e:v\in e}f(e)\leq1$, for every $v\in V(H).$ Furthermore, if equality holds for every $v\in V(H)$, then we call the fractional matching \textit{perfect}. Denote the maximum size of a fractional matching of $H$ by $\nu^*(H)=\max_f\Sigma_{e\in E(H)}f(e).$ Let $c_{k,d}$ be the minimum $d$-degree threshold for perfect fractional matchings in $k$-graphs, namely, for every $\varepsilon>0$ and sufficiently large $n\in\mathbb{N}$, every $n$-vertex $k$-graph $H$ with $\delta_d(H)\geq(c_{k,d}+\varepsilon)\binom{n-d}{k-d}$ contains a perfect fractional matching.

Aharoni, Holzman and Jiang \cite{2019Aharonifrac} studied the existence of rainbow fractional {matchings} in a $k$-graph system.
\begin{theorem}[\cite{2019Aharonifrac}]\label{fractional}
     Let $k\geq2$ be an integer, and let $n$ be a positive rational number. Assume that $\mathcal{H}=\{H_1,\ldots,H_{\lceil kn\rceil}\}$ is a $k$-graph system with $\nu^{*}(H_{i})\geq n$ for all $i\in [\lceil kn\rceil]$. Then $\mathcal{H}$ has a rainbow fractional matching of size $n.$   
\end{theorem}

When the $k$-graph is $k$-partite and $n$ is an integer, the number of graphs in $\mathcal{H}$ required in Theorem \ref{fractional} goes down from $kn$ to $kn-k+1.$ 
Theorem \ref{fractional} provided a fractional version of the corresponding problem concerning rainbow matchings, which was solved by Drisko \cite{1998Drisko} and by Aharoni and Berger \cite{2009AhaBerr-r} for  bipartite graphs (the rainbow matching problem remains open for general graphs as well as for 
$k$-partite hypergraphs when  $k>2$).

Alon, Frankl, Huang, R\"odl, Ruci\'nski and Sudakov \cite{2012Alonmatch}  established the relationship between perfect matching and perfect fractional matching by showing that every $n$-vertex $k$-graph $H$ with $\delta_d(H)\geq(\max\{c_{k,d},1/2\}+o(1))\binom{n-d}{k-d}$ has a perfect matching, and this condition is asymptotically best possible.  Cheng, Han, Wang and Wang \cite{2023chengspan} proved that the $d$-degree condition that ensures the existence of a perfect matching in a $k$-graph also guarantees the existence of a transversal perfect matching in a $k$-graph system for $d\in [k-1]$. 
\begin{theorem}[\cite{2023chengspan}]\label{thmwang}
     For every $\varepsilon > 0$ and integer $d\in [ k- 1]$,  there exists $n_{0}\in \mathbb{N}$,  such that the following holds for all integers $n\geq n_0$ and $n\in k\mathbb{N}$. Every $n$-vertex $k$-graph system $\mathcal{H}= \{ H_1, \ldots , H_{\frac nk}\}$  with  $\delta_d(\mathcal{H})\geq(\max\{c_{k,d},\frac12\}+\varepsilon)\binom{n-d}{k-d}$ contains a transversal perfect matching.  
\end{theorem}

Lang \cite{lang2023} gave a common generalization of the above theorem. Let $G_1$ and $G_2$ be $k$-graphs on the same vertex set $V$, and let $F$ be a $k$-graph. A \textit{color covering homomorphism} from $F$ to $(G_1,G_2)$ is a function $\phi:V(F)\rightarrow V$ with an edge $f\in E(F)$ such that $\phi(f)\in E(G_1)$ and $\phi$ is a homomorphism from $F'$ to $G_2$, where $F'$ is obtained from $F$ by deleting $f$ from its edge set. We define 
$\delta_d^{C_F}$ to be the smallest real number  $\delta$ such that for all $\alpha>0$ and $n$ large enough any two $k$-graphs $G_1, G_2$ on the same $n$ vertices with $\delta_d(G_1),\delta_d(G_2)\geq (\delta+\alpha)\binom{n-d}{k-d}$ admit a color covering homomorphism from $F$ to $(G_1,G_2)$.  Denote by $\delta_d^{F}$ the smallest real number $\delta$ such that for all $\alpha>0$ and {sufficiently large $n$}, every $n$-vertex $k$-graph $G$ with $\delta_d(G)\geq(\delta+\alpha)\binom{n-d}{k-d}$ contains a copy of $F$-factor.

\begin{theorem}[\cite{lang2023}]\label{Lang}
     Let $F$ be a $k$-graph with $r$ vertices and $t$ edges. For every $\varepsilon > 0$ and integer $d\in [k- 1]$, there exists $n_{0}\in \mathbb{N}$   such that the following holds for all integers $n\geq n_0$. Every $rn$-vertex $k$-graph system $\mathcal{H}= \{H_1, \ldots, H_{tn}\}$ with  $\delta_d(\mathcal{H})\geq\left(\max\left\{\delta_d^{C_F},\delta_d^{F}\right\}+\varepsilon\right)\binom{n-d}{k-d}$ contains a transversal $F$-factor.  
\end{theorem}

This is essentially a hypergraph version of the work of Montgomery, M\"{u}yesser and Pehova (see Theorem \ref{factor}). As an application of Theorem \ref{Lang}, Lang \cite{lang2023} deduced the minimum degree condition that guarantees the existence of transversal $K_4^3$-factors in a $3$-graph or transversal power of a Hamilton cycles in a $k$-graph, where $K_4^3$ is the complete $3$-graph on $4$ vertices. 

\section{Further results}
There are some additional results on transversal problems in graph systems whose details we will not discuss
in this survey. For convenience of the reader we give a concise  list of such
results together with relevant references.
\begin{itemize}
    \item 
    The rainbow Tur\'an problem can be viewed as a special case of Tur\'{a}n problems in graph systems.  Let $H$ be a graph. The \textit{rainbow Tur\'{a}n number} of $H$, denoted by $\mathrm{ex}^* (n, H)$, is defined to be the maximum number of edges in a properly edge-colored $n$-vertex graph $G$ with no rainbow copies of $H$. Notice that this properly edge-colored graph $G$ can be viewed as {a} system $\mathcal{G}=\{G_1,\ldots,G_{\chi(G)}\}$ of graphs where each $G_i$ is the graph consisting of edges with  color $i.$ Then {finding} a rainbow copy of $H$ in $G$ is equivalent to {finding} a transversal copy of $H$ inside  $\mathcal{G}$. The rainbow Tur\'{a}n number of some classical families of graphs were studied by Keevash, Mubayi, Sudakov and Verstra\"ete \cite{2007KeevashRTuran}. See more results in \cite{Bondy17974,2013DasC_2k,2017JohnstonRTu,2024KimRC}.

    \item In section 3, we {mainly} survey results on the minimum degree condition for transversal spanning structures in graph systems. In fact, some studies have considered the problem of embedding transversal non-spanning structures in a graph system. Li, Li and Li \cite{2022liluyi} studied the existence of transversal cliques in graph systems under Dirac-type conditions. Falgas-Ravry, Markstr\"om and R\"aty \cite{2024RavrymindT} studied the minimum degree version of Problem \ref{triangle}, and they found that the extremal behavior in the minimum degree condition differs strikingly from that seen in the density condition, with discrete jumps as opposed to continuous transitions. 

    \item The \textit{adjacency matrix} of $G$ is $A(G)=(a_{ij})$, where $a_{xy}=1$ if $xy\in E(G)$ and $0$ otherwise.  The largest eigenvalue of $A(G)$ is called the \textit{spectral radius} of $G$. The maximum spectral radius of graphs in a graph system that forbids transversal {perfect matchings, Hamilton paths or Hamilton cycles} was studied in \cite{2023GuoSpectral,he,zhangyuke}.
    \item One may consider the  problem of embedding a rainbow graph $H$ in a graph system $\mathcal{G}$ provided  that $H$ and all graphs in  $\mathcal{G}$ belong to the same class. For example, Aharoni, Briggs, Holzman and Jiang \cite{2021AharoniOddcyc}  proved that every family of $2\left\lceil n/2\right\rceil-1$ odd cycles of order $n$ contains a rainbow odd cycle. For more results, we refer the reader to \cite{2024DongEvencyc,2022GoorevitchTra,2023LiTreeStar}.
    {    \item  {Another classical and closely related concept is the independent transversal. Given a graph $G$ and a partition of its vertex set, an \textit{independent transversal} is an independent set of $G$ that contains one vertex from each partition class. It also concerns selecting one element from each of several sets while satisfying a global structural constraint.} 
   {  This is a very general notion, many combinatorial problems can be formulated by asking whether a certain graph with a certain vertex partition has an independent transversal. For more related results, see \cite{IT2024HaxellSiam,2025IT,IT2022Glock,IT2016Haxell,ITTopoMSur,IT2003Haxelltb,IT2006Szabo,IT1997Yuster}.}

}
\end{itemize}

\section{Conclusion}
Although we have made an effort to provide a systematic coverage of {recent 
results} in graph systems, 
it is inevitable that some  topics were left out of this survey due to constraints of space and time. Nevertheless, 
we would like to believe that we have presented enough results demonstrating how one can generalize known results in extremal graph theory  to {transversal versions}. 
Therefore we hope this survey will inspire further research in this area.

\section{Acknowledgement}
{We would like to express our gratitude to the anonymous reviewer for his/her valuable comments that greatly improved the presentation of this paper. We also thank Stijn Cambie, Yangyang Cheng,  Jie Hu and Richard Lang for helpful comments and valuable suggestions. }

\bibliographystyle{spmpsci}
\bibliography{survey}

\end{document}